\documentclass[12pt,reqno]{amsart}
\usepackage{amsfonts,amsmath,amssymb}

\allowdisplaybreaks
\advance\textheight by \topskip

\newmuskip\pFqskip
\mathchardef\pFcomma=\mathcode`, 

\newcommand*\pFq[5]{%
	\begingroup
	\begingroup\lccode`~=`,
	\lowercase{\endgroup\def~}{\pFcomma\mkern\pFqskip}%
	\mathcode`,=\string"8000
	{}_{#1}F_{#2}\biggl[\genfrac..{0pt}{}{#3}{#4}\Big|#5\biggr]%
	\endgroup
}

\usepackage{euscript}

\numberwithin{equation}{section}

\title{Two New Identities Involving the Catalan Numbers: A classical approach}
\author[H.~Prodinger]{Helmut Prodinger}
\address{H.~Prodinger\\Department of Mathematical Sciences, Mathematics Division\\ Stellenbosch University, Private Bag X1, 7602 Matieland, South Africa}
\email{hproding@sun.ac.za}

\thanks{}

\begin{document}
 \maketitle
 
 \section{Introduction}
 
 Mikic presents in \cite{Mikic} two identities involving Catalan numbers and provides combinatorial (bijective) proofs. Now, these days, when I hear  ``combinatorial identity'', I think immediately about Zeilberger's algorithm and hypergeometric series, and am curious, what they say. This is exactly the aim of the present note.
 Zeilberger's algorithm is described in his lovely book $A=B$~\cite{PWZ}.

 \section{The first identity}
 
 Catalan numbers are written as $C_k=\binom {2k}{k}\frac1{k+1}$. The first identity is: 
\begin{equation*}
f(n)=\sum_{k=0}^n(-1)^k\binom nk C_k\binom{2n-2k}{n-k}=\binom{n}{\lfloor n/2\rfloor}^2.
\end{equation*} 
 We ask Zeilberger's algorithm and get a recursion of second order:
 \begin{equation*}
-16(2n+1)(n-1)^2f(n-2)-4(2n^2-1)f(n-1)+f(n)(2n-1)(n+1)^2=0.
 \end{equation*}
 This is a bit surprising; of course, one could plug in the right hand side and prove the formula by induction.
 
 Let us write 
 \begin{equation*}
C_k=\binom{2k}{k}-\binom{2k}{k+1}
 \end{equation*}
 and treat the two terms separately.
\begin{equation*}
F(n)=\sum_{k=0}^n(-1)^k\binom nk \binom{2k}{k}\binom{2n-2k}{n-k}.
\end{equation*} 
 Now we get
 \begin{equation*}
F(n)n^2=-16(n-1)^2F(n-2), 
 \end{equation*}
 which can be iterated and leads to
 \begin{equation*}
F(2n)=\binom{2n}{n}^2,\quad F(2n+1)=0.
 \end{equation*}
 Similarly, with
 \begin{equation*}
 G(n)=\sum_{k=0}^n(-1)^k\binom nk \binom{2k}{k+1}\binom{2n-2k}{n-k}
 \end{equation*}
\begin{equation*}
 G(n)(n+1)^2=16n^2G(n-2)
 \end{equation*}
 and
 \begin{equation*}
G(2n+1)=-\binom{2n+1}{n}^2,\quad G(2n)=0.
 \end{equation*}
 So Zeilberger's algorithm strikes again; we only have to combine the two results for $F(n)$ and $G(n)$.
 
 But even without any special packages on sums my (old version of) Maple gives me for even $n$
 \begin{equation*}
F(n)=\frac{\Gamma^2\big(\frac{n+1}{2}\big)4^n}{\pi\Gamma^2\big(\frac{n}{2}+1\big)}
 \end{equation*}
 which is what it should be, using the duplication formula for the Gamma function. Similarly, we get for odd $n$
 \begin{equation*}
 G(n)=-\frac{\Gamma^2\big(\frac{n}{2}+1\big)4^n}{\pi\Gamma^2\big(\frac{n+3}{2}\big)}.
 \end{equation*}
 Now we switch to hypergeometric functions. Consider
\begin{equation*}
W_{k,l}(a,b,c):=\pFq{3}{2}{a,b,c}{\frac{1+a+b+k}{2},2c+l}{1}.
\end{equation*} 
As we find in Chu's paper \cite{Chu}, $W_{0,0}(a,b,c)$ can be evaluated by a classical formula of Watson in terms of 8 Gamma functions. This paper contains recursions about how to bring down $W_{k,l}(a,b,c)$ to the instance that can be evaluated. Here is the example we need; the notation of a Pochhammer symbol $(x)_n=x(x+1)\dots(x+n-1)$ is used.
 \begin{align*}
W_{0,1}(a,b,c)&=\sum_{j=0}^1(-1)^j\frac{(a)_j(b)_j}{(2c-1)_{2j}}
\frac{(c)_j(2c-1)_j}{(\frac{1+a+b}{2})_j(2c+1)_j}W_{0,0}(a+j,b+j,c+j)\\
&= W_{0,0}(a,b,c)-\frac{ab}{(1+a+b)(2c+1)}W_{0,0}(a+1,b+1,c+1).
 \end{align*}
 The sum in question is this:
 \begin{equation*}
 f(n)=\binom{2n}{n}\pFq{3}{2}{-n,-n,\tfrac12}{-n+\tfrac12,2}{1}.
 \end{equation*}
 So we see that we need $a=b=-n$, $c=\frac12$ and
 \begin{align*}
 W_{0,1}(-n,-n,\frac12)
 &= W_{0,0}(-n,-n,\frac12)+\frac{n^2}{2(2n-1)}W_{0,0}(-n+1,-n+1,\frac32).
 \end{align*}
Evaluating the appearing $W_{0,0}$ functions by Watson's formula leads to
  \begin{align*}
\binom{2n}{n} W_{0,1}(-n,-n,\frac12)
 &= \binom{2n}{n}W_{0,0}(-n,-n,\frac12)\\&+\binom{2n}{n}\frac{n^2}{2(2n-1)}W_{0,0}(-n+1,-n+1,\frac32)\\
 &=[\![n \text{ even} ]\!]\binom{n}{n/2}^2+[\![n \text{ odd} ]\!]\binom{n}{(n-1)/2}^2,
 \end{align*}
which is again the formula of interest. It is to be noted that modern computer algebra systems ``know'' some of these formul\ae, so that it is not necessary to type in Watson's formula.

\section{The second identity}
The formula is
\begin{equation*}
\sum_{k=0}^{2n}(-1)^k\binom{2n}{k}C_kC_{2n-k}=C_n\binom{2n}{n}.
\end{equation*}
So we set
\begin{equation*}
f(n):=\sum_{k=0}^n(-1)^k\binom nk C_k C_{n-k}
\end{equation*}
and ask Zeilberger:
\begin{equation*}
f(n)=\frac{16(n-1)^2}{n(n+2)}f(n-2),\qquad f(0)=C_{n}.
\end{equation*}
So we get nonzero values only for even $n$, and the announced formula follows by iteration/induction.

Now we move to hypergeometric functions. We need to evaluate
\begin{equation*}
f(n)=C_n\,\pFq{3}{2}{-n,-n-1,\tfrac12}{-n+\tfrac12,2}{1}.
\end{equation*}
This is a $W_{1,1}$ function in the notation of Chu~\cite{Chu}.

However, Maple can evaluate the function $\pFq{3}{2}{-n,-n-1,\tfrac12}{-n+\tfrac12,2}{1}$. It outputs an ugly version, but after using the reflection formula and the duplication formula for the Gamma function, it leads to 0 for odd $n$, and for $2n$ to
\begin{equation*}
\pFq{3}{2}{-2n,-2n-1,\tfrac12}{-2n+\tfrac12,2}{1}=\frac{(2n+1)!(2n)!^3}{(4n)!n!^3(n+1)!}.
\end{equation*}

\end{document}